\documentclass[pre]{imsart}

\RequirePackage{amsthm,amsmath}
\RequirePackage{natbib}
\RequirePackage[colorlinks,citecolor=blue,urlcolor=blue]{hyperref}
\RequirePackage{graphicx}
\RequirePackage{amsfonts}

\startlocaldefs
\numberwithin{equation}{section}
\theoremstyle{plain}
\newtheorem{theorem}{Theorem}
\newtheorem{proposition}[theorem]{Proposition}%
\newtheorem{corollary}[theorem]{Corollary}
\newtheorem{lemma}[theorem]{Lemma}
\theoremstyle{thmstylethree}%
\def\e{ \mathrm{e} }
\theoremstyle{remark}
\newtheorem{remark}{Remark}
\endlocaldefs

\begin{document}

\begin{frontmatter}
\title{Spectral clustering algorithm for the allometric extension model}
\runtitle{Spectral clustering for allometric extension}

\begin{aug}
\author{\fnms{Kohei} \snm{Kawamoto}\ead[label=e1]{kawamoto.kohei.532@s.kyushu-u.ac.jp}}

\address{Joint Graduate School of Mathematics for Innovation, Kyushu University\\
\printead{e1}}

\author{\fnms{Yuichi} \snm{Goto}\ead[label=e2]{yuichi.goto@math.kyushu-u.ac.jp}}
\and
\author{\fnms{Koji} \snm{Tsukuda}
\ead[label=e3]{tsukuda@math.kyushu-u.ac.jp}
}

\address{Faculty of Mathematics, Kyushu University\\
\printead{e2,e3}\\
}

\runauthor{K. Kawamoto et al.}

\end{aug}

\begin{abstract}
The spectral clustering algorithm is often used as a binary clustering method for unclassified data by applying the principal component analysis.
To study theoretical properties of the algorithm, the assumption of conditional homoscedasticity is often supposed in existing studies.
However, this assumption is restrictive and often unrealistic in practice.
Therefore, in this paper, we consider the allometric extension model, that is, the directions of the first eigenvectors of two covariance matrices and the direction of the difference of two mean vectors coincide, and we provide a non-asymptotic bound of the error probability of the spectral clustering algorithm for the allometric extension model.
As a byproduct of the result, we obtain the consistency of the clustering method in high-dimensional settings.
\end{abstract}


\end{frontmatter}

\section{{Introduction}}\label{sec:1}	

\subsection{Spectral clustering algorithm}\label{subsec:11}

Clustering unclassified data is one of typical problems in multivariate data analysis.
Several clustering methods have been intensively developed.
For example, the $k$-means clustering \citep{RefP81,RefP82}, the hierarchical clustering \citep{RefB14}, the method based on linear discriminant functions \citep{RefO78}, to mention only a few, are commonly used methods. 
For a review of clustering methods, see \cite{RefA17} and references therein. 
Of those methods, the spectral clustering algorithm, a clustering method based on principal component analysis, is quite popular because of its low computational load. 
It indeed has various applications such as community detection and graph partitioning. 
See \cite{RefL21} for references and further applications.
According to Section 4.7.1 of \cite{RefV18}, the method is described as follows.

Let the dimension of data and the sample size be $n$ and $m$, respectively, and $n$-dimensional (centered) data points are denoted by $\boldsymbol{x}_1, \ldots, \boldsymbol{x}_m$.
Calculate the realized sample covariance matrix
\[
\boldsymbol{S}_x = \frac{1}{m} \sum_{i=1}^m \boldsymbol{x}_i \boldsymbol{x}_i^\top
\]
and the unit-length eigenvector $\boldsymbol{v}$ corresponding to the largest eigenvalue of $\boldsymbol{S}_x$, where {$\top$ denotes the transpose of a vector.
To classify data points into two groups, we project data points onto the space spanned by $\boldsymbol{v}$  and classify them by the signs of their principal component scores.

To review theoretical properties of the spectral clustering algorithm, we define $\theta$ indicating group to which an individual belongs as a symmetric Bernoulli variable (Rademacher variable), that is to say, 
$P(\theta=1)=P(\theta=-1)=1/2$, and let $\boldsymbol{g}^{(1)}$ and  $\boldsymbol{g}^{(-1)}$ be $n$-dimensional centered normal variables with covariance matrices $\boldsymbol{\Sigma}_1$ and $\boldsymbol{\Sigma}_2$, respectively.
Define an $n$-dimensional random variable $\boldsymbol{X}$ as
\[ \boldsymbol{X} = \theta \boldsymbol{\mu} + \boldsymbol{g}^{(\theta)} ,\]
where $\boldsymbol{\mu}$ is a $n$-dimensional vector.
Note that $E[\boldsymbol{X}]= \boldsymbol{0}_n$.}
Let $\boldsymbol{X}_1,\ldots,\boldsymbol{X}_m$ be iid copies of $\boldsymbol{X}$, $\boldsymbol{S}_{m} = m^{-1} \sum_{i=1}^m \boldsymbol{X}_i \boldsymbol{X}_ i^\top$ the sample covariance matrix, and $\boldsymbol{\gamma}_1(\boldsymbol{S}_{m})$ the unit eigenvector  corresponding to the largest eigenvalue of $\boldsymbol{S}_{m}$.
Under the assumption that $\boldsymbol{\Sigma}_1 = \boldsymbol{\Sigma}_2 = \boldsymbol{I}_n$, where $\boldsymbol{I}_n$ is the the identity matrix of size $n$, \citet[Section 4.7]{RefV18} evaluated
\begin{equation} \label{verp}
P\left( \{\mbox{The number of misclassifications of $X_1,\ldots,X_m$}\} \leq \varepsilon m \right),
\end{equation}
where $\varepsilon$ is a positive constant satisfying appropriate conditions and 
``The number of misclassifications'' will be formulated later in \eqref{kerp}.
Moreover, under the assumption that $\boldsymbol{\Sigma}_1 = \boldsymbol{\Sigma}_2 = \boldsymbol{I}_n$ and $ \| \boldsymbol{\mu} \|_2 \geq C_{\rm gap} n/m $ for some large constant ${C_{\rm gap}}$, \cite{RefCZ18} derived an upper bound for the expectation of the misclustering rate that is a quantity of interest for clustering methods, where ``the misclustering rate''  will be explained in \eqref{misclustering}.
In this paper, we use the words ``misclassification'' and ``misclustering'' in the different senses.
\cite{RefDS13} proposed the clustering method based on the moment method for a mixture distribution of spherical normal distributions.
\cite{RefAFW22} proposed a more sophisticated method using the singular value decomposition of the Gram matrix calculated with its diagonal components replaced by zero and evaluated the misclustering rate of this method.
Recently, for a derivative algorithm of Lloyd's iterative procedure in a two-component mixture of normal distributions, \cite{RefN22} {derived conditions under which asymptotically correct results of the procedure are obtained}.

When $\boldsymbol{\Sigma}_1 = \boldsymbol{\Sigma}_2$, the first eigenvector of the covariance matrix $\boldsymbol{\Sigma}$ of the mixture distribution is different from that of $\boldsymbol{\Sigma}$ in general.
However, if $\boldsymbol{\Sigma}_1$ is spherical, the direction of the first eigenvector of $\boldsymbol{\Sigma}$ is parallel to $\boldsymbol{\mu}_1 - \boldsymbol{\mu}_2$, where $\boldsymbol{\mu}_1 = \boldsymbol{\mu}$ and $\boldsymbol{\mu}_2= -\boldsymbol{\mu}$.
This serves as the basis for the spectral clustering algorithm.
However, the assumption that $\boldsymbol{\Sigma}_1 = \boldsymbol{\Sigma}_2$ and $\boldsymbol{\Sigma}_1$ is spherical is sometimes restrictive and unlikely to hold in practice.
Note that Section 5 of \cite{RefAFW22} evaluated eigenvalues and eigenvectors in heteroscedastic situation for each individual.
Therefore, in this paper, we suppose a weaker assumption that $\boldsymbol{\mu}$, $\boldsymbol{\Sigma}_1$, and $\boldsymbol{\Sigma}_2$ follow the allometric extension relationship formulated in \cite{RefF97} and \cite{RefBFN99}, that is, the leading eigenvector of $\boldsymbol{\Sigma}_1$ is parallel to that of $\boldsymbol{\Sigma}_2$ and is also parallel to $\boldsymbol{\mu}_1-\boldsymbol{\mu}_2$. 
Then, we evaluate the misclassification probability of the spectral clustering algorithm. 
This results allows to show the consistency of the spectral clustering algorithm in a high-dimensional setting.

\subsection{{Notation}} 
In this subsection, we introduce some notations used in this paper.
For vectors $\boldsymbol{a}$ and $\boldsymbol{b}$ which have the same dimension, $\langle \boldsymbol{a}, \boldsymbol{b} \rangle = {\boldsymbol{a}^{\top} \boldsymbol{b} }$ denotes the inner product of  $\boldsymbol{a}$ and $\boldsymbol{b}$.
For a vector $\boldsymbol{a}$, let $\| \boldsymbol{a} \|_2 = \sqrt{\langle \boldsymbol{a}, \boldsymbol{a} \rangle}$ be the $\ell^2$-norm of $\boldsymbol{a}$.
For a square matrix $\boldsymbol{A}$, let $\lambda_k(\boldsymbol{A})$ be the $k$-th largest eigenvalue of $\boldsymbol{A}$, $\boldsymbol{\gamma}_k(\boldsymbol{A})$ the eigenvector corresponding to $\lambda_k(\boldsymbol{A})$, and $\| \boldsymbol{A} \|_{\mathrm{op}}$ the operator norm of $\boldsymbol{A}$.
For a positive integer $n$, $\mathcal{N}_n(\boldsymbol{\mu},\boldsymbol{\Sigma})$ denotes the $n$-dimensional normal distribution with a mean vector $\boldsymbol{\mu}$ and a covariance matrix $\boldsymbol{\Sigma}$.
Moreover, in general, $(\boldsymbol{\mu},\boldsymbol{\Sigma})$ denotes the distribution with a mean vector $\boldsymbol{\mu}$ and a covariance matrix $\boldsymbol{\Sigma}$.
For a positive definite covariance matrix $\boldsymbol{\Sigma}$, $\boldsymbol{\Sigma}^{1/2}$ denotes the symmetric matrix satisfying $\boldsymbol{\Sigma}^{1/2} \boldsymbol{\Sigma}^{1/2} = \boldsymbol{\Sigma}$, and $\boldsymbol{\Sigma}^{- 1/2}$ denotes the inverse of  $\boldsymbol{\Sigma}^{1/2}$.
For a random variable $X$, $\| X \|_{\psi_2}$ and $\| X \|_{L^2}$ denote the sub-gaussian norm and the $L^2$ norm of $X$, respectively.

\subsection{Organization of {the} paper}
In section~\ref{sec:2}, the allometric extension model is introduced and some properties of this model are derived. 
Section~\ref{sec:3} provides the main results of this paper.
Theoretical results are proven in Section~\ref{sec:4}.
Some concluding remarks are presented in section~\ref{sec:5}.

\section{Allometric {e}xtension {m}odel}\label{sec:2}
In this section, we introduce the allometric extension model. 
Let $n$ be a positive integer, and let $\boldsymbol{\mu}_1$ and $\boldsymbol{\mu}_2$ be $n$-dimensional vectors satisfying $\boldsymbol{\mu}_1 \neq \boldsymbol{\mu}_2$. 
Also, let $\boldsymbol{\Sigma}_1$ and $\boldsymbol{\Sigma}_2$ be $n \times n$ positive definite symmetric matrices such that $\lambda_1(\boldsymbol{\Sigma}_i) > \lambda_2(\boldsymbol{\Sigma}_i)$ for $i=1,2$. 
We define the allometric extension relationship between two distributions $(\boldsymbol{\mu}_1, \boldsymbol{\Sigma}_1)$ and $(\boldsymbol{\mu}_2, \boldsymbol{\Sigma}_2)$ as follows: there exists a $\beta \in \mathbb{R}$ such that
\begin{align}\label{aem}
\boldsymbol{\gamma}_1(\boldsymbol{\Sigma}_1) = \boldsymbol{\gamma}_1(\boldsymbol{\Sigma}_2) = \beta (\boldsymbol{\mu}_1 - \boldsymbol{\mu}_2) ,
\end{align}
where the sign of $\boldsymbol{\gamma}_1(\boldsymbol{\Sigma}_2)$ is suitably chosen. 
Throughout the discussion, we assume that the two distribution are $n$-dimensional normal distributions.

The allometric extension model formulated in \cite{RefF97} and \cite{RefBFN99} are used as a model to express a typical relationship between two or more biotic groups. 
For this model, \cite{RefBFN99} proposed a test procedure for determining whether two groups are in the allometric extension relationship and analyzed carapace size data of turtles of different sexes as discussed in the paper by Jolicoeur and Mosimann; for this data, we refer to Table 1.4 in \cite{RefF97}.
Moreover, \cite{RefTM23} proposed a test procedure for the allometric extension relationship when observations are high-dimensional. 
Furthermore, properties of mixtures of two or more distributions that form the allometric extension relationship are also discussed in \citet[section 8.7]{RefF97}, \cite{RefKHF08} and \cite{RefMK14}. 
Note that \cite{RefT07} mentioned the case of applying $k$-means algorithm to the allometric extension model.
We focus on the spectral clustering algorithm, which is computationally more efficient than several traditional clustering algorithms due to its dimensionality reduction techniques.

Let us derive some properties associated with the mixture distribution of two normal distributions forming the allometric extension relationship.
Let $n$ be a positive integer, and let $f_{1}(\cdot)$ and $f_{2}(\cdot)$ be the probability density functions of $\mathcal{N}_n(\boldsymbol{\mu}_1,\boldsymbol{\Sigma}_1)$ and $\mathcal{N}_n(\boldsymbol{\mu}_2,\boldsymbol{\Sigma}_2)$, respectively. 
We assume that {$(\boldsymbol{\mu}_1, \boldsymbol{\Sigma}_1)$ and $(\boldsymbol{\mu}_2, \boldsymbol{\Sigma}_2)$ satisfy  \eqref{aem}} and that the $n$-dimensional random variable $\boldsymbol{X}$ follows a mixture distribution whose probability density function is given by
\[ f_X(\boldsymbol{x}) = \pi_1 f_{1}(\boldsymbol{x}) + \pi_2 f_{2}(\boldsymbol{x}) \quad (\boldsymbol{x} \in \mathbb{R}^n), \]
where $\pi_1$ and $\pi_2$ are positive values satisfying $\pi_1 + \pi_2 = 1$.
In this case, the following properties about the covariance matrix $\boldsymbol{\Sigma} = V[\boldsymbol{X}]$ hold.

\begin{proposition}\label{propae}
Under the conditions stated above,
\begin{align}
\mbox{(i)} & \quad \boldsymbol{\gamma}_1(\boldsymbol{\Sigma}) = \boldsymbol{\gamma}_1(\boldsymbol{\Sigma}_1),  \label{AE1} \\
\mbox{(ii)} & \quad \lambda_1(\boldsymbol{\Sigma} ) = \pi_1 \lambda_1(\boldsymbol{\Sigma}_1)+ \pi_2 \lambda_1(\boldsymbol{\Sigma}_2)+ \pi_1 \pi_2 \| \boldsymbol{\mu}_1 - \boldsymbol{\mu}_2 \|_2^2, \label{AE2} \\
\mbox{(iii)} & \quad \lambda_2(\boldsymbol{\Sigma}) \leq \pi_1 \lambda_2(\boldsymbol{\Sigma}_1) + \pi_2 \lambda_2(\boldsymbol{\Sigma}_2) , \label{AE3}
\end{align}
where the sign of $\boldsymbol{\gamma}_1(\boldsymbol{\Sigma})$ is appropriately chosen in \eqref{AE1}.
\end{proposition}

\begin{remark}
The results \eqref{AE1} and \eqref{AE2} are given in Lemma 8.7.1 of \cite{RefF97}, but \eqref{AE3} gives a better evaluation than the corresponding result in \cite{RefF97}.
Therefore, we include the proof of Proposition~\ref{propae}, although the proof is similar.
\end{remark}

\begin{remark}
In this study, we consider the situation that $\lambda_1(\boldsymbol{\Sigma}_i) > \lambda_2(\boldsymbol{\Sigma}_i)$ for $i=1,2$. 
Therefore, we stress that the case where $\boldsymbol{\Sigma}_1 = \boldsymbol{\Sigma}_2 = \boldsymbol{I}_n$ is not included in our setting.
\end{remark}

\section{{Spectral clustering algorithm for the allometric extension model}}\label{sec:3}
{In this section, we derive an upper bound on the misclassification probability of the spectral clustering algorithm for the allometric extension model.} 
Let $\theta$ be a symmetric Bernoulli variable, $n$ a positive integer, $\boldsymbol{\mu}$ an $n$-dimensional vector, and $\boldsymbol{\Sigma}_1, \boldsymbol{\Sigma}_2$ $n \times n$ positive definite symmetric matrices satisfying $\lambda_1(\boldsymbol{\Sigma}_i) > \lambda_2(\boldsymbol{\Sigma}_i)$ for $i=1,2$. 
Suppose that
\[ \boldsymbol{\gamma}_1(\boldsymbol{\Sigma}_1) =  \boldsymbol{\gamma}_1(\boldsymbol{\Sigma}_2) = 2\beta \boldsymbol{\mu} \]
for some $\beta > 0$.
We consider an $n$-dimensional random variable $\boldsymbol{X}$ defined as
\begin{align*}
\boldsymbol{X} = \theta \boldsymbol{\mu} + \boldsymbol{g}^{(\theta)},
\end{align*}
{where $\boldsymbol{g}^{(1)} \sim \mathcal{N}_n(\boldsymbol{0}_n, \boldsymbol{\Sigma}_1)$ and $\boldsymbol{g}^{(-1)} \sim \mathcal{N}_n(\boldsymbol{0}_n, \boldsymbol{\Sigma}_2)$ that are independent of $\theta$}. Then, the probability density function $f_X(\cdot)$ of $\boldsymbol{X}$ is given by
\[ 
f_{X}(\boldsymbol{x}) = \frac{1}{2} f_{1}(\boldsymbol{x}) + \frac{1}{2} f_{2}(\boldsymbol{x}) \quad (\boldsymbol{x} \in \mathbb{R}^n),
\]
where $f_{1}(\cdot)$ and $f_{2}(\cdot)$ are the probability density functions of $\mathcal{N}_n(\boldsymbol{\mu},\boldsymbol{\Sigma}_1)$ and $\mathcal{N}_n(-\boldsymbol{\mu},\boldsymbol{\Sigma}_2)$, respectively.
As shown in the following proposition, $\boldsymbol{X}$ is an $\mathbb{R}^n$-valued sub-gaussian random variable.

\begin{proposition}\label{kex}
There exist a constant $K \geq 1$ such that
\begin{equation}\label{kthi}
\| \langle \boldsymbol{X}, \boldsymbol{x} \rangle \|_{\psi_2} \leq K \| \langle \boldsymbol{X} , \boldsymbol{x} \rangle \|_{L^2}
\end{equation}
for any $\boldsymbol{x} \in \mathbb{R}^n$.
In particular, \eqref{kthi} holds for any $\boldsymbol{x} \in \mathbb{R}^n$ when $K = \sqrt{{32}/(4 - \e) } (= 4.9966\cdots)$.
\end{proposition}

For a positive integer $m$, we observe random variables $\boldsymbol{X}_1,\ldots,\boldsymbol{X}_m$ following 
\begin{equation*}
\boldsymbol{X}_i = \theta_i \boldsymbol{\mu} + \boldsymbol{g}_i^{(\theta_i)} \quad (i=1,\ldots,m),
\end{equation*}
where $\theta_1,\ldots,\theta_m$ are iid copies of $\theta$ and $\boldsymbol{g}^{(t)}_1,\ldots,\boldsymbol{g}^{(t)}_m$ are iid copies of $\boldsymbol{g}^{(t)}$ for $t =-1,1$.
Note that $\boldsymbol{X}_1,\ldots,\boldsymbol{X}_m$ are iid copies of $\boldsymbol{X}$.
We presume that there are two groups $\mathcal{N}_n(\boldsymbol{\mu},\boldsymbol{\Sigma}_1)$ and $\mathcal{N}_n(-\boldsymbol{\mu},\boldsymbol{\Sigma}_2)$ forming the allometric extension relationship and $\theta_i$ indicates the group to which an individual $\boldsymbol{X}_i$ belongs for $i=1,\ldots,m$.
The sample covariance matrix $\boldsymbol{S}_{m} $ is defined as
\begin{equation*}
\boldsymbol{S}_{m} = \frac{1}{m} \sum_{i=1}^m \boldsymbol{X}_i \boldsymbol{X}_i^\top.
\end{equation*}
Here we note that $E[\boldsymbol{X}_i] = \boldsymbol{0}_n$ for $i=1,\ldots,m$.
As an estimator of 
\begin{equation}\label{aes}
\boldsymbol{\Sigma}= E[\boldsymbol{X} \boldsymbol{X}^\top] 
= \boldsymbol{\mu} \boldsymbol{\mu}^\top + \frac{1}{2} \boldsymbol{\Sigma}_1 + \frac{1}{2} \boldsymbol{\Sigma}_2, 
\end{equation}
the estimation error $\boldsymbol{S}_m - \boldsymbol{\Sigma}$ of $\boldsymbol{S}_m$ can be evaluated in the following proposition, which will be used to prove our main result.

\begin{proposition}\label{evs}
Let $K (\geq 1)$ be a constant satisfying \eqref{kthi} for any $\boldsymbol{x} \in \mathbb{R}^n$.
For any $u \geq 0$,
\begin{align*}
&P\left(\| \boldsymbol{ S }_m - \boldsymbol{  \Sigma } \|_{\mathrm{op}} \leq C K^2 \left( \sqrt{ \frac{n+u}{m} }  + \frac{n+u}{m} \right) \left(\frac{\lambda_1(\boldsymbol{\Sigma}_1) + \lambda_1(\boldsymbol{\Sigma}_2)}{2} + \| \boldsymbol{\mu} \|_2^2 \right)
\right)\\
&\geq 1- 2\e^{-u},
\end{align*}
where $C$ is some positive absolute constant. 
\end{proposition}

\begin{remark}
If we use the result of Exercise 4.7.3 of \cite{RefV18}, Proposition~\ref{evs} immediately follows from Propositions~\ref{propae} and \ref{kex}.
The proof of Proposition~\ref{evs} is included just to be sure.
\end{remark}

\begin{remark}
It also holds that
\[
E[ \| \boldsymbol{ S }_m - \boldsymbol{  \Sigma } \|_{\mathrm{op}} ]
\leq C K^2 \left( \sqrt{\frac{n}{m}} + \frac{n}{m} \right) \left( \frac{\lambda_1(\boldsymbol{\Sigma}_1) + \lambda_1(\boldsymbol{\Sigma}_2)}{2} + \| \boldsymbol{\mu} \|_2^2 \right) ,
\]
which is a direct consequence of Theorem 4.7.1 of \cite{RefV18} and Propositions~\ref{propae} and \ref{kex}.
\end{remark}

If we regard binary clustering as a binary classification problem for unlabeled data, the main objective is to classify individuals in a random sample into correct groups, where each individual belongs to one of the two groups.
Let us assume that $\langle \boldsymbol{\gamma}_1(\boldsymbol{S}_{m}), \boldsymbol{\gamma}_1(\boldsymbol{\Sigma}) \rangle >0$.
In our problem setting, the spectral clustering algorithm classifies $\boldsymbol{X}_1,\ldots,\boldsymbol{X}_m$ into two clusters by the signs of $\langle \boldsymbol{\gamma}_1(\boldsymbol{S}_{m}), \boldsymbol{X}_i \rangle $ $(i=1,\ldots,m)$.
The misclassification probability of $\boldsymbol{X}_i$ can be expressed as
\[ P(\theta_i \langle \boldsymbol{\gamma}_1(\boldsymbol{S}_{m}), \boldsymbol{X}_i \rangle < 0) \]
for $i=1,\ldots,m$.
As it will be remarked later, evaluating misclassification probability enables us to evaluate the misclustering rate; see Remark~\ref{remmisclust}.
The following theorem provides a non-asymptotic upper bound of this misclassification probability.
 
\begin{theorem}\label{mthm}
Let 
$C$, $K$, $K_g$, and $c$ be positive absolute constants which {are independent of } $m$, $n$, $\boldsymbol{\mu}$, $\boldsymbol{\Sigma}_1$, and $\boldsymbol{\Sigma}_2$, and let $c_1 = 1+K_g^2/\sqrt{c}$.
{Suppose that} $m$, $n$, $\boldsymbol{\mu}$, $\boldsymbol{\Sigma}_1$, and $\boldsymbol{\Sigma}_2$ satisfy
\begin{align}
&\sqrt{2} C K^2 \left ( \sqrt{\frac{2n}{m}} + \frac{2n}{m} \right) \left( \frac{\lambda_1(\boldsymbol{\Sigma}_1) +  \lambda_1(\boldsymbol{\Sigma}_2)}{ \| \boldsymbol{\mu} \|_2^2} + 2 \right)
\nonumber\\ 
 &\leq \frac{\alpha \| \boldsymbol{\mu}\|_2}{c_1 \sqrt{n \max_{j=1,2}\{\lambda_1(\boldsymbol{\Sigma}_j)\}} + \| \boldsymbol{\mu} \|_2}
 \label{tha}
\end{align}
for some $\alpha \in (0,1)$. 
Then it holds that
\begin{align*}
P(\theta_i \langle \boldsymbol{\gamma}_1(\boldsymbol{S}_{m}), \boldsymbol{X}_i \rangle <0)  &\leq \Phi\left( \frac{-(1-\alpha)\|\boldsymbol{\mu}\|_2}{\sqrt{\max_{j=1,2} \{ \lambda_1(\boldsymbol{\Sigma}_j) \} }} \right) + 6 \e^{-n}
\end{align*}
for $i=1,\ldots,m$, where $\Phi(\cdot)$ is the distribution function of the standard normal distribution.
\end{theorem}

\begin{remark}
Let us explain the constants $C$, $K$, $K_g$, and $c$ in Theorem~\ref{mthm}.
The constant $K$ is given in Proposition~\ref{kex}, and $C$ in Proposition~\ref{evs}.
Moreover, the constant $K_g$ is the sub-gaussian norm of a standard normal variable; in particular, $K_g = \sqrt{8/3}$.
Letting $\boldsymbol{g} \sim \mathcal{N}_n (\boldsymbol{0}_n, \boldsymbol{I}_n)$, we have
\begin{equation}\label{v33}
P( \bigl| \| \boldsymbol{g}\|_2  -\sqrt{n} \bigr| \geq t) \leq 2 \exp\left( - \frac{ct^2}{K_g^4} \right)
\end{equation}
for all $t\geq 0$; see, e.g.,  Equation (3.3) in \citet[p.40]{RefV18}.
The constant $c$ in Theorem~\ref{mthm} appears in the inequality \eqref{v33}.
\end{remark}

Theorem \ref{mthm} provides a non-asymptotic lower bound for the probability of the misclassification rate announced in Section~\ref{sec:1}.
In our setting, \eqref{verp} is formulated as
\begin{align} 
&P\left( \sum_{i=1}^m 1\{ \theta_i \langle \boldsymbol{\gamma}_1(\boldsymbol{S}_{m}), \boldsymbol{X}_i \rangle < 0 \} \leq \varepsilon m \right)
\nonumber\\
&= 1 - P\left( \sum_{i=1}^m 1\{ \theta_i \langle \boldsymbol{\gamma}_1(\boldsymbol{S}_{m}), \boldsymbol{X}_i \rangle < 0 \} > \varepsilon m \right).
\label{kerp}
\end{align}
By applying the Markov inequality, we have
\begin{align*}
&P\left( \sum_{i=1}^m 1\{ \theta_i \langle \boldsymbol{\gamma}_1(\boldsymbol{S}_{m}), \boldsymbol{X}_i \rangle < 0 \} > \varepsilon m \right) \\
&\leq \frac{1}{\varepsilon m} \sum_{i=1}^m  E\left[1\left\{ \theta_i \langle \boldsymbol{\gamma}_1(\boldsymbol{S}_{m}), \boldsymbol{X}_i \rangle < 0 \right\} \right]
= \frac{1}{\varepsilon} P\left( \theta_1 \langle \boldsymbol{\gamma}_1(\boldsymbol{S}_{m}), \boldsymbol{X}_1 \rangle < 0 \right),
\end{align*}
because $(\theta_1 \langle \boldsymbol{\gamma}_1(\boldsymbol{S}_{m}), \boldsymbol{X}_1 \rangle, \ldots, \theta_m \langle \boldsymbol{\gamma}_1(\boldsymbol{S}_{m}), \boldsymbol{X}_m \rangle )$ is exchangeable.
Thus, we obtain the following corollary to Theorem~\ref{mthm}.

\begin{corollary}\label{co0}
Consider constants $C$, $K$, and $c_1$ in Theorem~\ref{mthm}. 
{Suppose that} $m$, $n$, $\boldsymbol{\mu}$, $\boldsymbol{\Sigma}_1$, and $\boldsymbol{\Sigma}_2$ satisfy \eqref{tha} for some $\alpha \in (0,1)$.  Then it holds that
\begin{equation*} 
P\left( \sum_{i=1}^m 1\{ \theta_i \langle \boldsymbol{\gamma}_1(\boldsymbol{S}_{m}), \boldsymbol{X}_i \rangle < 0 \} \leq \varepsilon m \right)
\geq 1 - \frac{1}{\varepsilon}\Phi\left( \frac{-(1-\alpha)\|\boldsymbol{\mu}\|_2}{\sqrt{\max_{j=1,2} \{\lambda_1(\boldsymbol{\Sigma}_j)\}}} \right) - \frac{6}{\varepsilon} \e^{-n} .
\end{equation*}
\end{corollary}

\begin{remark}\label{remmisclust}
In our setting, the misclustering rate stated in Section~\ref{sec:1} is expressed as
\begin{equation}\label{misclustering}
\frac{1}{m} \min \left\{ \sum_{i=1}^m 1\{ \theta_i \langle \boldsymbol{\gamma}_1(\boldsymbol{S}_m), \boldsymbol{X}_i \rangle < 0\}, \sum_{i=1}^m 1\{ \theta_i \langle \boldsymbol{\gamma}_1(\boldsymbol{S}_m), \boldsymbol{X}_i \rangle > 0\}  \right\}.
\end{equation}
It holds that
\begin{align*}
& P\left( \frac{1}{m} \min \left\{ \sum_{i=1}^m 1\{ \theta_i \langle \boldsymbol{\gamma}_1(\boldsymbol{S}_m), \boldsymbol{X}_i \rangle < 0\}, \sum_{i=1}^m 1\{ \theta_i \langle \boldsymbol{\gamma}_1(\boldsymbol{S}_m), \boldsymbol{X}_i \rangle > 0\}  \right\} \leq \varepsilon \right) \\
& \geq P\left( \sum_{i=1}^m 1\{ \theta_i \langle \boldsymbol{\gamma}_1(\boldsymbol{S}_m), \boldsymbol{X}_i \rangle < 0\} \leq  \varepsilon m \right).
\end{align*}
We can evaluate the right-hand side by using Corollary~\ref{co0}. 
In addition, the expectation of \eqref{misclustering} may be similarly evaluated.
Note that \eqref{misclustering} is called ``misclassification rate'' in \cite{RefCZ18} and \cite{RefAFW22}.
\end{remark}

We define a signal-to-noise ratio $\eta$ for our problem as 
\[ \eta = \frac{\|\boldsymbol{\mu}\|_2^2}{\max_{j=1,2} \{ \lambda_1(\boldsymbol{\Sigma}_j)\}}, \]
which plays an important role to evaluate performances.
The numerator $\|\boldsymbol{\mu}\|_2^2$ and the denominator $\max_{j=1,2}\{\lambda_1(\boldsymbol{\Sigma}_j)\}$ correspond to strength of signals and noises, respectively. 
If noises are large compared to signals, then signals are hidden by noises, which makes detection difficult.
In this sense, $\eta$ can be interpreted as the difficulty of classification when the allometric extension model is considered. 
If $\eta$ is small (large), it is difficult (easy) to classify individuals in a random sample into the correct groups. 
It is easy to see that \eqref{tha} is fulfilled when $m$, $n$, $\boldsymbol{\mu}$, $\boldsymbol{\Sigma}_1$, and $\boldsymbol{\Sigma}_2$ satisfy
\begin{equation}\label{thas}
2^{ 3/2 } C K^2 \left ( \sqrt{\frac{2n}{m}} + \frac{2n}{m} \right) \left( \frac{1}{\eta} + 1 \right) 
 \leq \frac{\alpha }{c_1 \sqrt{n/\eta} + 1}.
\end{equation}
Hence, if $n/\eta = O(1)$, then \eqref{tha} holds for sufficiently large $m$ compared to $n$.

\begin{remark}
\cite{RefAFW22} and \cite{RefN22} discuss this kind of topic in a similar context.
\end{remark}

\begin{remark}
By the Mills inequality, Theorem~\ref{mthm} leads to
\begin{align} 
&P(\theta_i \langle \boldsymbol{\gamma}_1(\boldsymbol{S}_{m}), \boldsymbol{X}_i \rangle \leq 0) 
\nonumber\\
&\leq
\frac{1}{\sqrt{2\pi (1-\alpha)^{2}\eta}} \exp\left(- \frac{(1-\alpha)^2 \eta}{2} \right) + 6 \e^{-n}
\quad (i=1,\ldots,m).
\label{millsa}
\end{align}
The inequality \eqref{millsa} shows that the misclassification probability for an individual decays exponentially with the signal-to-noise ratio $\eta$ (multiplied by a constant) and the dimension $n$.
\end{remark}

Finally, we consider the probability of the event
\begin{equation*}
\left\{ \bigcap_{i=1}^m \bigl\{ \theta_i \langle \boldsymbol{\gamma}_1(\boldsymbol{S}_{m}), \boldsymbol{X}_i \rangle > 0 \bigr\} \right\} \cup \left\{ \bigcap_{i=1}^m \bigl\{ \theta_i \langle \boldsymbol{\gamma}_1(\boldsymbol{S}_{m}), \boldsymbol{X}_i \rangle <0 \bigr\} \right\}
\end{equation*}
indicating all individuals $\boldsymbol{X}_1,\ldots,\boldsymbol{X}_m$ are clustered correctly, that is to say, the misclustering rate equals zero.
Theorem~\ref{mthm} implies the consistency of clustering in the sense of \eqref{consc} under a high-dimensional regime.

\begin{corollary}\label{col}
As $m,n \to \infty$ with
\begin{equation}\label{ARh}
 \frac{n}{m} \to 0, \quad  \frac{\log{m}}{n} \to 0 ,
\end{equation}
if $n/\eta = O(1)$, then
\begin{equation}\label{consc}
P \left( \left\{ \bigcap_{i=1}^m \bigl\{ \theta_i \langle \boldsymbol{\gamma}_1(\boldsymbol{S}_{m}), \boldsymbol{X}_i \rangle > 0 \bigr\} \right\} \cup \left\{ \bigcap_{i=1}^m \bigl\{ \theta_i \langle \boldsymbol{\gamma}_1(\boldsymbol{S}_{m}), \boldsymbol{X}_i \rangle <0 \bigr\} \right\}  \right) \to 1.
\end{equation}
\end{corollary}

\section{Proofs}\label{sec:4}

\subsection{Proof of Proposition \ref{propae}}

For simplicity, denote $\boldsymbol{\gamma}_1 (\boldsymbol{\Sigma }_1) = \boldsymbol{\gamma}_1 ( \boldsymbol{\Sigma}_2)$ by $\boldsymbol{\gamma}_1 $, where the sign of $\boldsymbol{\gamma}_1 (\boldsymbol{\Sigma}_2)$ is appropriately chosen.
Let us construct a random variable $Y$ taking values in $\{1,2\}$ such that $P(Y = i) = \pi_i$ and $\boldsymbol{X} | \{Y=i\} \sim \mathcal{N}_n (\boldsymbol{\mu}_i, \boldsymbol{\Sigma}_i)$ conditionally for $i=1,2$.
Then, we have  $E[ \boldsymbol{X} \arrowvert Y = i] = \mu_i$ and $V[\boldsymbol{X} \arrowvert Y = i] = \boldsymbol{  \Sigma }_i$ for $i=1,2$.
These formulae give
\begin{equation*}
E[V[\boldsymbol{X}\arrowvert Y]] = \sum_{i=1}^2 \pi_i \boldsymbol{ \Sigma }_i, \quad
V[ E[ \boldsymbol{X} \arrowvert Y]]= \sum_{i=1}^2 \pi_i \left( \boldsymbol{\mu}_i - \bar{\boldsymbol{\mu}} \right) \left( \boldsymbol{\mu}_i - \bar{\boldsymbol{\mu}} \right)^\top,
\end{equation*}
where
\[ \bar{\boldsymbol{\mu}} = \sum_{i=1}^2 \pi_i\boldsymbol{\mu}_i = \boldsymbol{\mu}_1 - \frac{\pi_2}{\beta} \boldsymbol{\gamma}_1 
=\boldsymbol{\mu}_2 + \frac{\pi_1}{\beta} \boldsymbol{\gamma}_1 .\]
It follows that
\begin{align*}
&\boldsymbol{  \Sigma } 
= E[V[\boldsymbol{X}\arrowvert Y]] + V[ E[ \boldsymbol{X} \arrowvert Y]] 
= \sum_{i=1}^2 \pi_i \boldsymbol{ \Sigma }_i + \sum_{i=1}^2 \pi_i \left( \boldsymbol{\mu}_i - \bar{\boldsymbol{\mu}} \right) \left( \boldsymbol{\mu}_i - \bar{\boldsymbol{\mu}} \right)^\top 
\\
&=  \sum_{i=1}^2 \pi_i \boldsymbol{ \Sigma }_i +  \frac{\pi_1 \pi_2}{\beta^2} \boldsymbol{\gamma}_1  \boldsymbol{\gamma}_1 ^\top ,
\end{align*}
because
\begin{align*}
&\sum_{i=1}^2 \pi_i \left( \boldsymbol{\mu}_i - \bar{\boldsymbol{\mu}} \right) \left( \boldsymbol{\mu}_i - \bar{\boldsymbol{\mu}} \right)^\top 
= \pi_1 \left(\frac{\pi_2}{\beta}\boldsymbol{\gamma}_1 \right) \left(\frac{\pi_2}{\beta}\boldsymbol{\gamma}_1 \right)^\top
+ \pi_2 \left(-\frac{\pi_1}{\beta}\boldsymbol{\gamma}_1 \right) \left( -\frac{\pi_1}{\beta}\boldsymbol{\gamma}_1 \right)^\top \\*
&=\frac{\pi_1\pi_2^2}{\beta^2} \boldsymbol{\gamma}_1  \boldsymbol{\gamma}_1 ^\top +
 \frac{\pi_1^2 \pi_2}{\beta^2}\boldsymbol{\gamma}_1 \boldsymbol{\gamma}_1 ^\top
= \frac{\pi_1 \pi_2}{\beta^2} (\pi_1 + \pi_2)\boldsymbol{\gamma}_1 \boldsymbol{\gamma}_1 ^\top 
= \frac{\pi_1 \pi_2}{\beta^2} \boldsymbol{\gamma}_1  \boldsymbol{\gamma}_1 ^\top,
\end{align*}
which is due to the definition of the allometric extension model.
From
\begin{equation*}
\boldsymbol{  \Sigma } \boldsymbol{\gamma}_1  
= \sum_{i=1}^2 \pi_i \boldsymbol{ \Sigma }_i \boldsymbol{\gamma}_1  + \frac{\pi_1 \pi_2}{\beta^2} \boldsymbol{\gamma}_1  \boldsymbol{\gamma}_1 ^\top \boldsymbol{\gamma}_1 
= \left( \sum_{i=1}^2 \pi_i \lambda_1 ( \boldsymbol{ \Sigma }_i )  +   \frac{ \pi_1  \pi_2}{\beta^2} \right) \boldsymbol{\gamma}_1 ,
\end{equation*}
we see that $\boldsymbol{\gamma}_1 $ is an eigenvector of $\boldsymbol{  \Sigma }$ corresponding to the eigenvalue 
\begin{equation*}
\sum_{i=1}^2 \pi_i       \lambda_1 ( \boldsymbol{ \Sigma }_i )  + \frac{\pi_1  \pi_2}{\beta^2}.
\end{equation*}
Let $\boldsymbol{ \xi }$ be another unit-length eigenvector of $\boldsymbol{\Sigma}$ orthogonal to $\boldsymbol{\gamma}_1$.
As $\boldsymbol{ \xi }$ and $\boldsymbol{\gamma}_1 $ are orthogonal, it holds that
\begin{equation*}
\boldsymbol{ \xi }^\top \boldsymbol{  \Sigma } \boldsymbol{ \xi }
= \sum_{i=1}^2 {\pi_i} \boldsymbol{ \xi }^\top \boldsymbol{ \Sigma }_i \boldsymbol{ \xi } + \frac{ \pi_1 \pi_2 }{\beta^2} \boldsymbol{ \xi }^\top \boldsymbol{\gamma}_1  \boldsymbol{\gamma}_1 ^\top \boldsymbol{ \xi }
= \sum_{i=1}^2{ \pi_i } \boldsymbol{ \xi }^\top \boldsymbol{ \Sigma }_i \boldsymbol{ \xi } .
\end{equation*}
By using a spectral decomposition
\begin{equation*}
\boldsymbol{\Sigma }_i  = \lambda_1( \boldsymbol{\Sigma }_i ) \boldsymbol{\gamma}_1  \boldsymbol{\gamma}_1 ^\top + \sum_{j=2}^n
\lambda_j( \boldsymbol{\Sigma}_i ) \boldsymbol{\gamma}_j ( \boldsymbol{\Sigma}_i ) \boldsymbol{\gamma}_j( \boldsymbol{\Sigma}_i ) ^\top
\quad (i = 1,2),
\end{equation*}
we have
\begin{align*}
&\boldsymbol{ \xi }^\top \boldsymbol{  \Sigma }_i \boldsymbol{ \xi } 
= \boldsymbol{ \xi }^\top \left( \sum_{j=2}^n \lambda_j( \boldsymbol{  \Sigma }_i ) \boldsymbol{\gamma}_j( \boldsymbol{  \Sigma }_i ) \boldsymbol{\gamma}_j( \boldsymbol{  \Sigma }_i )^\top \right)  \boldsymbol{ \xi } 
=  \sum_{j=2}^n \lambda_j ( \boldsymbol{  \Sigma }_i ) \langle \boldsymbol{ \xi } , \boldsymbol{\gamma}_j( \boldsymbol{  \Sigma }_i ) \rangle ^2 \\
& \leq \lambda_2( \boldsymbol{  \Sigma }_i ) \sum_{j=2}^n  \langle \boldsymbol{ \xi } , \boldsymbol{\gamma}_j( \boldsymbol{  \Sigma }_i ) \rangle ^2 
=  \lambda_2( \boldsymbol{  \Sigma }_i )
\end{align*}
for $i=1,2$.
Consequently, we deduce that
\[
\boldsymbol{ \xi }^\top \boldsymbol{  \Sigma } \boldsymbol{ \xi } \leq \sum_{i=1}^2  \pi_i \lambda_2(\boldsymbol{  \Sigma }_i)
< \sum_{i=1}^2 \pi_i \lambda_1 ( \boldsymbol{ \Sigma }_i )  +   \frac{ \pi_1  \pi_2}{\beta^2} 
\]
for any unit-length eigenvector $\boldsymbol{\xi}$ of $\boldsymbol{\Sigma}$ orthogonal to $\boldsymbol{\gamma}_1$, hence 
\begin{equation*}
\sum_{i=1}^2 \pi_i \lambda_1 ( \boldsymbol{ \Sigma }_i )  + \frac{\pi_1  \pi_2}{\beta^2}
\end{equation*}
is the unique largest eigenvalue $\lambda_1(\boldsymbol{\Sigma})$ of $\boldsymbol{\Sigma}$ and 
\begin{equation*}
\lambda_2(\boldsymbol{  \Sigma }) \leq \sum_{i=1}^2 \pi_i \lambda_2(\boldsymbol{  \Sigma }_i) 
.
\end{equation*}
This completes the proof.
\qed

\subsection{Proof of Proposition~\ref{kex}}

When $\boldsymbol{x} = \boldsymbol{0}_n$, \eqref{kthi} holds for any $K (\geq 1)$.
Hereafter, we consider $\boldsymbol{x} \neq \boldsymbol{0}_n$.

We first observe that
\begin{equation*}
\| \langle \boldsymbol{X}, \boldsymbol{x} \rangle \|_{L^2}^2 
= E [\langle \boldsymbol{X}, \boldsymbol{x} \rangle^2] 
= E [\boldsymbol{x}^{\top} \boldsymbol{X}\boldsymbol{X}^{\top} \boldsymbol{x}] 
= \boldsymbol{x}^{\top} E[\boldsymbol{X}\boldsymbol{X}^{\top}] \boldsymbol{x} = \langle \boldsymbol{\Sigma} \boldsymbol{x}, \boldsymbol{x} \rangle.
\end{equation*}
Next, we evaluate $\| \langle \boldsymbol{X}, \boldsymbol{x} \rangle \|_{\psi_2}$.
It holds that
\begin{align}
& E \biggl[\exp \left( \frac{ { \langle \boldsymbol{X}, \boldsymbol{x} \rangle }^2 }{ t^2 } \right) \biggr] 
=  E \biggl[ \exp \left( \frac{ { \langle \theta \boldsymbol{\mu} + \boldsymbol{g}^{(\theta)}, \boldsymbol{x} \rangle }^2 }{ t^2 } \right) \biggr] \nonumber \\
&= E \biggl[ E\biggl[ \exp \left( \frac{ { \langle \theta \boldsymbol{\mu} + \boldsymbol{g}^{(\theta)}, \boldsymbol{x} \rangle }^2 }{ t^2 } \right) \biggl| \theta \biggr] \biggr] \nonumber \\
&= \frac{1}{2} E \biggl[ \exp{ \left( \frac{ { \langle \boldsymbol{\mu} + \boldsymbol{g}^{(1)}, \boldsymbol{x} \rangle }^2 }{ t^2 } \right)} \biggr] 
+ \frac{1}{2} E \biggl[ \exp{ \left( \frac{ { \langle -\boldsymbol{\mu} + \boldsymbol{g}^{(-1)}, \boldsymbol{x} \rangle }^2 }{ t^2 } \right)} \biggr] ,
\label{kthp1}
\end{align}
where $t$ will be specified later.
As for the first term on the right-hand side of \eqref{kthp1}, it follows from $\langle \boldsymbol{\Sigma}_1^{- 1/2} \boldsymbol{g}^{(1)} , { \boldsymbol{\Sigma}_1^{ 1/2}\boldsymbol{x} }/{\| \boldsymbol{\Sigma}_1^{ 1/2}\boldsymbol{x} \|_2} \rangle \sim \mathcal{N}(0,1)$ that
\begin{align*}
& \frac{1}{2} E  \left[ \exp{ \left( \frac{  \langle \boldsymbol{\mu} + \boldsymbol{g}^{(1)}, \boldsymbol{x} \rangle ^2 }{ t^2 } \right) }  \right]  \\
&\leq \frac{1}{2} E \left[ \exp { \left( \frac{  {2\langle \boldsymbol{\mu} , \boldsymbol{x} \rangle }^2  + 2 {\langle  \boldsymbol{g}^{(1)}, \boldsymbol{x} \rangle }^2 } { t^2 } \right)}
\right] \\
&= \frac{1}{2} \exp{ \left( \frac{ 2\langle \boldsymbol{\mu} , \boldsymbol{x} \rangle^2 }{ t^2 } \right) }  E \left[   \exp{ \left( \frac{ 2\langle \boldsymbol{\Sigma}_1^{- 1/2} \boldsymbol{g}^{(1)} , \boldsymbol{\Sigma}_1^{ 1/2}\boldsymbol{x} \rangle^2 }{ t^2 } \right) }      \right] \\
&= \frac{1}{2} \exp{ \left( \frac{ 2\langle \boldsymbol{\mu} , \boldsymbol{x} \rangle^2 }{ t^2 } \right) }  E \left[   \exp{ \left( \frac{ 2 \| \boldsymbol{\Sigma}_1^{ 1/2}\boldsymbol{x} \|_2^2 \left\langle \boldsymbol{\Sigma}_1^{- 1/2} \boldsymbol{g}^{(1)} , 
{ \boldsymbol{\Sigma}_1^{ 1/2}\boldsymbol{x} }/{\| \boldsymbol{\Sigma}_1^{ 1/2}\boldsymbol{x} \|_2} \right\rangle^2 }{ t^2 } \right) }      \right] \\
&= \frac{1}{2} \exp{ \left( \frac{ 2\langle \boldsymbol{\mu} , \boldsymbol{x} \rangle^2 }{ t^2 } \right) } E \left[ \exp \left( \frac{ Z^2 }{ { t^2 }/{( 2\| \boldsymbol{\Sigma}_1^{ 1/2}\boldsymbol{x} \|_2^2 )} } \right) \right], 
\end{align*}
where $Z$ is a standard normal random variable.
The expectation on the right-hand side is finite when 
$t^2 > 4\| \boldsymbol{\Sigma}_1^{ 1/2}\boldsymbol{x} \|_2^2$ 
and is given by
\begin{equation*}
E \left[ \exp \left( \frac{ Z^2 }{ { t^2 }/{( 2\| \boldsymbol{\Sigma}_1^{ 1/2}\boldsymbol{x} \|_2^2 )} } \right) \right] 
= \frac{1}{ \sqrt{ 1 - {4\| \boldsymbol{\Sigma}_1^{ 1/2}\boldsymbol{x} \|_2^2 }/{ t^2 }} },
\end{equation*}
which yields that
\begin{equation*}
\frac{1}{2} E \biggl[ \exp{ \left( \frac{  \langle \boldsymbol{\mu} + \boldsymbol{g}^{(1)}, \boldsymbol{x} \rangle ^2 }{ t^2 } \right) }  \biggr] \leq \frac{1}{2} \exp{ \left( \frac{ 2\langle \boldsymbol{\mu} , \boldsymbol{x} \rangle^2 }{ t^2 } \right)}  \frac{1}{ \sqrt{ 1 - {4\| \boldsymbol{\Sigma}_1^{ 1/2}\boldsymbol{x} \|_2^2 }/{ t^2 }} }.
\end{equation*}
Similarly, as for the second term on the right-hand side of \eqref{kthp1}, when $t^2 > 4\| \boldsymbol{\Sigma}_2^{ 1/2}\boldsymbol{x} \|_2^2$, it holds that
\begin{equation*}
\frac{1}{2} E  \biggl[ \exp{ \left( \frac{  \langle -\boldsymbol{\mu} + \boldsymbol{g}^{(-1)}, \boldsymbol{x} \rangle ^2 }{ t^2 } \right) }  \biggr] \leq \frac{1}{2} \exp{ \left( \frac{ 2\langle \boldsymbol{\mu} , \boldsymbol{x} \rangle^2 }{ t^2 } \right) }  \frac{1}{ \sqrt{ 1 - {4\| \boldsymbol{\Sigma}_2^{ 1/2}\boldsymbol{x} \|_2^2 }/{ t^2 }} }.
\end{equation*}
Letting $ M = \max\{ \| \boldsymbol{\Sigma}_1^{1/2} \boldsymbol{x} \|_2 , \| \boldsymbol{\Sigma}_2^{1/2} \boldsymbol{x} \|_2\} $, we have
\begin{equation*}
E \biggl[\exp \left( \frac{ { \langle \boldsymbol{X}, \boldsymbol{x} \rangle }^2 }{ t^2 } \right) \biggr] \leq   \exp{ \left( \frac{ 2\langle \boldsymbol{\mu} , \boldsymbol{x} \rangle^2 }{ t^2 } \right) } \frac{1}{ \sqrt{ 1 - {4  M^2 }/{ t^2 }} } 
\end{equation*}
when $t^2> 4M^2$.
Hence, from the definition of the sub-gaussian norm, it suffices to show that there exists a constant $K$ satisfying this inequality with $ t = K \sqrt{ \langle \boldsymbol{\Sigma} \boldsymbol{x}, \boldsymbol{x} \rangle }$ for all $\boldsymbol{x} (\neq \boldsymbol{0}_n)$.
It follows from \eqref{aes} that 
\begin{equation*}
\langle \boldsymbol{\Sigma} \boldsymbol{x} , \boldsymbol{x} \rangle 
= \langle \boldsymbol{\mu} , \boldsymbol{x} \rangle^2+ \frac{1}{2} \| \boldsymbol{\Sigma}_1^{1/2} \boldsymbol{x} \|_2^2 + \frac{1}{2} \| \boldsymbol{\Sigma}_2^{1/2} \boldsymbol{x} \|_2^2 
\geq \langle \boldsymbol{\mu} , \boldsymbol{x} \rangle^2.
\end{equation*}
Moreover, plugging $t = K \sqrt{ \langle \boldsymbol{\Sigma} \boldsymbol{x}, \boldsymbol{x} \rangle }$ into $\exp ( { 2\langle \boldsymbol{\mu} , \boldsymbol{x} \rangle^2 }/{ t^2 })$, we have
\begin{equation*}
\exp \left( \frac{ 2\langle \boldsymbol{\mu} , \boldsymbol{x} \rangle^2 }{ t^2 } \right)  
= \exp \left( \frac{ 2\langle \boldsymbol{\mu} , \boldsymbol{x} \rangle^2 }{ K^2 \langle \boldsymbol{\Sigma} \boldsymbol{x}, \boldsymbol{x} \rangle  } \right)
\leq 
\exp \left( \frac{ 2\langle \boldsymbol{\mu} , \boldsymbol{x} \rangle^2 }{ K^2 \langle \boldsymbol{\mu} , \boldsymbol{x} \rangle^2 } \right) 
= \exp \left( \frac{2}{K^2} \right).
\end{equation*}
Furthermore, if $K \geq 2$, then $\exp({2}/{K^2}) \leq \sqrt{\e}$.
Therefore, what is left is to show that there exists a constant $K \geq 2$ such that
\begin{equation}\label{kthp2}
\frac{1}{ \sqrt{ 1 - {4 M^2 }/(K^2 \langle \boldsymbol{\Sigma} \boldsymbol{x}, \boldsymbol{x} \rangle )} }
\leq \frac{2}{ \sqrt{\e} } \end{equation}
for all $\boldsymbol{x} (\neq 0)$.
The inequality \eqref{kthp2} is equivalent to
\[ K^2 \geq \frac{16}{\e} \left( \frac{4}{\e} -  1 \right)^{-1} \frac{M^2}{ \langle \boldsymbol{\Sigma} \boldsymbol{x} , \boldsymbol{x} \rangle } .\]
From
\begin{equation*}
\frac{M^2}{ \langle \boldsymbol{\Sigma} \boldsymbol{x} , \boldsymbol{x} \rangle } 
= \frac{M^2}{ \langle \boldsymbol{\mu} , \boldsymbol{x} \rangle^2+ \| \boldsymbol{\Sigma}_1^{1/2} \boldsymbol{x} \|_2^2 / 2 +  \| \boldsymbol{\Sigma}_2^{1/2} \boldsymbol{x} \|_2^2 / 2 }
\leq \frac{2M^2}{  \| \boldsymbol{\Sigma}_1^{1/2} \boldsymbol{x} \|_2^2 +  \| \boldsymbol{\Sigma}_2^{1/2} \boldsymbol{x} \|_2^2 } 
\leq 2,
\end{equation*}
it follows that 
\[ K =  \sqrt{\frac{32}{\e} \left( \frac{4}{\e} -  1 \right)^{-1}} = \sqrt{\frac{32}{4 - \e} } \]
satisfies \eqref{kthp2} for all $\boldsymbol{x} (\neq \boldsymbol{0}_n)$.
\qed

\subsection{Proof of Proposition~\ref{evs}}

Fix $u \geq 0$.
Let
\begin{equation*}
\boldsymbol{Z}_i = \boldsymbol{\Sigma}^{-1/2} \boldsymbol{X}_i \quad (i=1,\ldots,m), \quad \mbox{and} \quad
\boldsymbol{R} = \frac{1}{m} \sum_{i=1}^m \boldsymbol{Z}_i \boldsymbol{Z}_i^\top - \boldsymbol{I}_n.
\end{equation*}
It follows from Equation (4.25) of \citet[p.94]{RefV18} that
\begin{equation*}
\| \boldsymbol{S}_m - \boldsymbol{\Sigma} \|_{\mathrm{op}} \leq \| \boldsymbol{R} \|_{\mathrm{op}} \|\boldsymbol{\Sigma} \|_{\mathrm{op}}.
\end{equation*}
By using Proposition~\ref{kex} and Equation (4.22) of \citet[p.91]{RefV18} with $t=\sqrt{u}$, we have
\begin{equation*}
P\left( \| \boldsymbol{R} \|_{\mathrm{op}} \leq K^2 \max\{\delta,\delta^2\} \right) \geq
1 - 2\e^{-u},
\end{equation*}
where $\delta = \tilde{C}(\sqrt{n} + \sqrt{u})/\sqrt{m}$.
Note that $\tilde{C}$ is the absolute constant $C$ in Equation (4.22) of \citet{RefV18}.
We see that
\begin{equation*}
 \delta = \tilde{C} \biggl( \frac{ \sqrt{n} + \sqrt{u}}{ \sqrt{m}} \biggr)  \leq  \tilde{C} \sqrt{ \frac{2(n + u)}{m} } .
\end{equation*}
Letting $C = \max \{ \sqrt{2} \tilde{C}, 2 \tilde{C}^2 \}$, we have
\begin{equation*}
K^2 \max\{\delta,\delta^2\}
\leq  K^2(\delta + \delta^2)\leq C K^2 \biggl( \sqrt{ \frac{n+u}{m}} + \frac{n+u}{m}\biggr).
\end{equation*}
Consequently, it follows that
\begin{align*}
& P \left( \|  \boldsymbol{S}_m - \boldsymbol{\Sigma} \|_{\mathrm{op}} \leq C K^2 \biggl( \sqrt{ \frac{n+u}{m}} + \frac{n+u}{m} \biggr) \|\boldsymbol{\Sigma} \|_{\mathrm{op}}
\right) \\
&\geq 
P \left( \| \boldsymbol{R} \|_{\mathrm{op}} \leq C K^2 \biggl( \sqrt{ \frac{n+u}{m}} + \frac{n+u}{m} \biggr)  \right) \\
& \geq 
P \left( \| \boldsymbol{R} \|_{\mathrm{op}} \leq K^2 \max\{\delta,\delta^2\}  \right)
\geq 1 - 2 \e^{-u}.
\end{align*}
Finally, \eqref{AE2} implies  $\lambda_1(\boldsymbol{\Sigma} ) = \left(\lambda_1(\boldsymbol{\Sigma}_1)+  \lambda_1(\boldsymbol{\Sigma}_2)\right)/2+  \| \boldsymbol{\mu} \|_2^2$.
This completes the proof.
\qed

\subsection{Proof of Theorem \ref{mthm}}
We see that
\begin{align*}
& P\biggl( \theta_i \langle \boldsymbol{\gamma}_1(\boldsymbol{S}_{m}) , \boldsymbol{X}_i\rangle < 0 \biggr) \\*
&= P \biggl( \langle \boldsymbol{\gamma}_1(\boldsymbol{S}_{m}) , \boldsymbol{X}_i \rangle<0 \biggl\arrowvert \theta_i = 1 \biggr) P\biggl( \theta_i = 1 \biggr) \\*
& \quad + P \biggl( \langle \boldsymbol{\gamma}_1(\boldsymbol{S}_{m}) , \boldsymbol{X}_i  \rangle>0  \biggl\arrowvert \theta_i = -1 \biggr) P\biggl(\theta_i = -1 \biggr) \\*
&= \frac{1}{2} P \biggl( \langle \boldsymbol{\gamma}_1(\boldsymbol{S}_{m}) , \boldsymbol{X}_i\rangle<0 \biggl\arrowvert \theta_i = 1 \biggr) + 
\frac{1}{2}  P \biggl( \langle \boldsymbol{\gamma}_1(\boldsymbol{S}_{m}) , \boldsymbol{X}_i \rangle>0 \biggl\arrowvert \theta_i = -1 \biggr) \\*
&\leq \frac{1}{2} \left( \Phi\left( \frac{-(1-\alpha) \| \boldsymbol{\mu} \|_2 }{\sqrt {\lambda_1(\boldsymbol{\Sigma}_1)}} \right) + \Phi\left( \frac{-(1-\alpha) \| \boldsymbol{\mu} \|_2}{\sqrt{\lambda_1(\boldsymbol{\Sigma}_2)}} \right) \right) + 6 \e^{-n} 
 \\*
& \leq  \Phi\left( \frac{-( 1 - \alpha ) \| \boldsymbol{\mu} \|_2}{\sqrt{\max_{j=1,2} \{ \lambda_1(\boldsymbol{\Sigma}_j) \}}} \right) + 6 \e^{-n} ,
\end{align*}
the second last inequality being a consequence of Lemma~\ref{lemk} presented below.
\qed

\begin{lemma}\label{lemk}
Consider constants $C$, $K$, and $c_1$ in Theorem~\ref{mthm}. 
If $m$, $n$, $\boldsymbol{\mu}$, $\boldsymbol{\Sigma}_1$, and $\boldsymbol{\Sigma}_2$ satisfy
\[
\sqrt{2} C K^2 \left ( \sqrt{\frac{2n}{m}} + \frac{2n}{m} \right) \left( \frac{\lambda_1(\boldsymbol{\Sigma}_1) +  \lambda_1(\boldsymbol{\Sigma}_2)}{ \| \boldsymbol{\mu} \|_2^2} + 2 \right) 
 \leq \frac{\alpha \| \boldsymbol{\mu} \|_2}{c_1 \sqrt{n \lambda_1(\boldsymbol{\Sigma}_1)} + \| \boldsymbol{\mu} \|_2}
\]
for some $\alpha \in (0,1)$, then {it holds that}
\begin{equation}\label{lem2a}
P(\langle \boldsymbol{\gamma}_1(\boldsymbol{S}_{m}), \boldsymbol{X}_i \rangle < 0 | \theta_i = 1)
\leq \Phi\left( \frac{-(1-\alpha)\|\boldsymbol{\mu}\|_2}{\sqrt{\lambda_1(\boldsymbol{\Sigma}_1)}} \right) + 6 \e^{-n}
\end{equation}
for $i=1,\ldots,m$.
If $m$, $n$, $\boldsymbol{\mu}$, $\boldsymbol{\Sigma}_1$, and $\boldsymbol{\Sigma}_2$ satisfy
\[
\sqrt{2} C K^2 \left ( \sqrt{\frac{2n}{m}} + \frac{2n}{m} \right) \left( \frac{\lambda_1(\boldsymbol{\Sigma}_1) +  \lambda_1(\boldsymbol{\Sigma}_2)}{ \| \boldsymbol{\mu} \|_2^2} + 2 \right) 
 \leq \frac{\alpha \| \boldsymbol{\mu}\|_2}{c_1 \sqrt{n \lambda_1(\boldsymbol{\Sigma}_2)} + \| \boldsymbol{\mu} \|_2}
\]
for some $\alpha \in (0,1)$, then {it holds that}
\begin{equation}\label{lem2b}
P(\langle \boldsymbol{\gamma}_1(\boldsymbol{S}_{m}), \boldsymbol{X}_i \rangle > 0 | \theta_i = - 1)
\leq \Phi\left( \frac{-(1-\alpha)\|\boldsymbol{\mu}\|_2}{\sqrt{\lambda_1(\boldsymbol{\Sigma}_2)}} \right) + 6 \e^{-n}
\end{equation}
for $i=1,\ldots,m$.
\end{lemma}

This lemma provides an evaluation of the conditional misclassification probability, which is the key to prove Theorem~\ref{mthm}.
In the next subsection, let us prove Lemma~\ref{lemk}.

\subsection{Proof of Lemma~\ref{lemk}}\label{subsec;lem3}
We shall only prove \eqref{lem2a}, because the proof of \eqref{lem2b} is similar.

Fix $i \in \{1,\ldots,m \}$.
Let us denote the event $\{  \theta_i = 1\}$ by $\mathrm{A}_i$ for simplicity.
We have
\begin{align}
&P \biggl(\langle \boldsymbol{\gamma}_1(\boldsymbol{S}_{m}) , \boldsymbol{X}_i \rangle<0 \biggl\arrowvert \mathrm{A}_i\biggr) \nonumber \\*
&= P \biggl(\langle \boldsymbol{\gamma}_1(\boldsymbol{S}_{m}) - \boldsymbol{\gamma}_1(\boldsymbol{\Sigma}_1) + \boldsymbol{\gamma}_1(\boldsymbol{\Sigma}_1) , \boldsymbol{X}_i - \boldsymbol{\mu} + \boldsymbol{\mu}\rangle<0 \biggl\arrowvert \mathrm{A}_i \biggr) \nonumber \\*
&= P\biggl(  \langle \boldsymbol{\gamma}_1(\boldsymbol{S}_{m}) - \boldsymbol{\gamma}_1 (\boldsymbol{\Sigma}_1), \boldsymbol{X}_i - \boldsymbol{\mu} \rangle + \langle \boldsymbol{\gamma}_1(\boldsymbol{S}_{m}) -\boldsymbol{\gamma}_1(\boldsymbol{\Sigma}_1),\boldsymbol{\mu}\rangle +\langle          \boldsymbol{\gamma}_1(\boldsymbol{\Sigma}_1) ,\boldsymbol{X}_i - \boldsymbol{\mu}\rangle \nonumber \\
& \qquad +\langle\boldsymbol{\gamma}_1(\boldsymbol{\Sigma}_1),\boldsymbol{\mu}\rangle <0 \biggl\arrowvert \mathrm{A}_i\biggr) \nonumber \\*
&\leq  
P\biggl(\langle\boldsymbol{\Sigma}_1^{1/2}(\boldsymbol{\gamma}_1(\boldsymbol{S}_{m})- \boldsymbol{\gamma}_1(\boldsymbol{\Sigma}_1)),\boldsymbol{\Sigma}_1^{-1/2}(\boldsymbol{X}_i - \boldsymbol{\mu})\rangle - \|\boldsymbol{\gamma}_1(\boldsymbol{S}_{m}) - \boldsymbol{\gamma}_1(\boldsymbol{\Sigma}_1)\|_2\cdot \| \boldsymbol{\mu} \|_2 \nonumber \\*
&\qquad +  \|\boldsymbol{\Sigma}_1^{1/2}\boldsymbol{\gamma}_1(\boldsymbol{\Sigma}_1)\|_2 \biggl\langle \frac{\boldsymbol{\Sigma}_1^{1/2}\boldsymbol{\gamma}_1(\boldsymbol{\Sigma}_1)}{\| \boldsymbol{\Sigma}_1^{1/2}\boldsymbol{\gamma}_1(\boldsymbol{\Sigma}_1)\|_2 },\boldsymbol{\Sigma}_1^{-1/2}(\boldsymbol{X}_i - \boldsymbol{\mu})\biggr\rangle
+ \| \boldsymbol{\mu}\|_2 < 0\biggl\arrowvert \mathrm{A}_i \biggr) \nonumber \\*
&\leq  
P\biggl(
-\|\boldsymbol{\Sigma}_1^{1/2} \|_{\mathrm{op}} \|\boldsymbol{\gamma}_1(\boldsymbol{S}_{m}) - \boldsymbol{\gamma}_1 (\boldsymbol{\Sigma}_1)\|_2  \|\boldsymbol{\Sigma}_1^{- 1/2}(\boldsymbol{X}_i - \boldsymbol{\mu})\|_2
 \nonumber \\*
&\qquad 
 - \|\boldsymbol{\gamma}_1(\boldsymbol{S}_{m}) - \boldsymbol{\gamma}_1(\boldsymbol{\Sigma}_1)\|_2\cdot \| \boldsymbol{\mu} \|_2 \nonumber \\
 & \qquad
 +  \|\boldsymbol{\Sigma}_1^{1/2}\boldsymbol{\gamma}_1(\boldsymbol{\Sigma}_1)\|_2 \biggl\langle \frac{\boldsymbol{\Sigma}_1^{1/2}\boldsymbol{\gamma}_1(\boldsymbol{\Sigma}_1)}{\| \boldsymbol{\Sigma}_1^{1/2}\boldsymbol{\gamma}_1(\boldsymbol{\Sigma}_1)\|_2 },\boldsymbol{\Sigma}_1^{-1/2}(\boldsymbol{X}_i - \boldsymbol{\mu})\biggr\rangle
+ \| \boldsymbol{\mu}\|_2 < 0\biggl\arrowvert \mathrm{A}_i \biggr) 
, \label{lem2p1}
\end{align} 
where
$\langle \boldsymbol{\gamma}_1(\boldsymbol{\Sigma}_1),\boldsymbol{\mu} \rangle = \| \boldsymbol{\mu} \|_2$
and
\begin{align*}
&\langle \boldsymbol{\gamma}_1(\boldsymbol{S}_{m}) - \boldsymbol{\gamma}_1 (\boldsymbol{\Sigma}_1), \boldsymbol{X}_i - \boldsymbol{\mu} \rangle  
= \langle \boldsymbol{\Sigma}_1^{1/2} (\boldsymbol{\gamma}_1(\boldsymbol{S}_{m}) - \boldsymbol{\gamma}_1 (\boldsymbol{\Sigma}_1)), \boldsymbol{\Sigma}_1^{-1/2}(\boldsymbol{X}_i - \boldsymbol{\mu} )\rangle  \\
&\geq -
   \|\boldsymbol{\Sigma}_1^{1/2} (\boldsymbol{\gamma}_1(\boldsymbol{S}_{m}) - \boldsymbol{\gamma}_1 (\boldsymbol{\Sigma}_1))\|_2 \|\boldsymbol{\Sigma}_1^{- 1/2}(\boldsymbol{X}_i - \boldsymbol{\mu} )\|_2   \\
&\geq -\|\boldsymbol{\Sigma}_1^{1/2} \|_{\mathrm{op}} \|\boldsymbol{\gamma}_1(\boldsymbol{S}_{m}) - \boldsymbol{\gamma}_1 (\boldsymbol{\Sigma}_1)\|_2  \|\boldsymbol{\Sigma}_1^{- 1/2}(\boldsymbol{X}_i - \boldsymbol{\mu})\|_2 
 \end{align*}
are used.
Let $\boldsymbol{g}_i  = \boldsymbol{\Sigma}_1^{-1/2}(\boldsymbol{X}_i - \boldsymbol{\mu}) $ and 
\[ Z = \frac{1}{\sqrt{\lambda_1(\boldsymbol{\Sigma}_1)}} \|\boldsymbol{\Sigma}_1^{1/2}\boldsymbol{\gamma}_1(\boldsymbol{\Sigma}_1)\|_2 \biggl\langle \frac{\boldsymbol{\Sigma}_1^{1/2}\boldsymbol{\gamma}_1(\boldsymbol{\Sigma}_1)}{\|\boldsymbol{\Sigma}_1^{1/2}\boldsymbol{\gamma}_1(\boldsymbol{\Sigma}_1)\|_2 },\boldsymbol{g}_i \biggr\rangle.\]
Then, it can be shown that 
$\sqrt{\lambda_1(\boldsymbol{\Sigma}_1)} Z | \mathrm{A}_i \sim \mathcal{N}(0, \lambda_1(\boldsymbol{\Sigma}_1))$
 conditionally because 
$\boldsymbol{g}_i | \mathrm{A}_i \sim \mathcal{N}_{n}(0, \boldsymbol{I}_n)$
conditionally and the normal distribution enjoys the property of reproductivity.
The right-hand side of \eqref{lem2p1} is equal to
\begin{align*}
&P\Biggl( - \| \boldsymbol{\gamma}_1(\boldsymbol{S}_{m}) - \boldsymbol{\gamma}_1(\boldsymbol{\Sigma}_1)\|_2 \left(\sqrt{\lambda_1(\boldsymbol{\Sigma}_1)}\| \boldsymbol{g}_i\|_2 + \|\boldsymbol{\mu} \|_2\right) + \sqrt{\lambda_1(\boldsymbol{\Sigma}_1)} \cdot Z + \|\boldsymbol{\mu} \|_2 \\
&\qquad
<0\biggl\arrowvert \mathrm{A}_i \Biggr)\\
&=P\Biggl(Z<  \| \boldsymbol{\gamma}_1(\boldsymbol{S}_{m}) - \boldsymbol{\gamma}_1(\boldsymbol{\Sigma}_1)\|_2 \left(\| \boldsymbol{g}_i \|_2 + \frac{\| \boldsymbol{\mu} \|_2}{\sqrt{\lambda_1(\boldsymbol{\Sigma}_1)}} \right) - \frac{\| \boldsymbol{\mu} \|_2}{\sqrt{\lambda_1(\boldsymbol{\Sigma}_1)}}\biggl\arrowvert \mathrm{A}_i \Biggr).
\end{align*}
Let $\delta$ be a constant satisfying
\begin{equation*}
\sqrt{2} C K^2 \biggl( \sqrt{ \frac{2n}{m} } + \frac{2n}{m} \biggr) \biggl( \frac{\lambda_1(\boldsymbol{\Sigma}_1)}{\| \boldsymbol{\mu} \|_2^2} + \frac{\lambda_1(\boldsymbol{\Sigma}_2)}{\| \boldsymbol{\mu} \|_2^2}  + 2 \biggr) 
\leq  \delta \leq
\frac{\alpha     \| \boldsymbol{\mu} \|_2 /\sqrt{\lambda_1(\boldsymbol{\Sigma}_1)}}{ c_1 \sqrt{n} +   \| \boldsymbol{\mu} \|_2 / \sqrt{\lambda_1(\boldsymbol{\Sigma}_1)}}.
\end{equation*}
Then, we see that 
\begin{align}
& P\Biggl(Z<  \| \boldsymbol{\gamma}_1(\boldsymbol{S}_{m}) - \boldsymbol{\gamma}_1(\boldsymbol{\Sigma}_1)\|_2 \left(\| \boldsymbol{g}_i \|_2 + \frac{\| \boldsymbol{\mu} \|_2}{\sqrt{\lambda_1(\boldsymbol{\Sigma}_1)}} \right) - \frac{\| \boldsymbol{\mu} \|_2}{\sqrt{\lambda_1(\boldsymbol{\Sigma}_1)}} \biggl\arrowvert \mathrm{A}_i \Biggr) \nonumber \\
&= P\Biggl( \Biggl\{Z<  \|\boldsymbol{\gamma}_1(\boldsymbol{S}_{m}) - \boldsymbol{\gamma}_1(\boldsymbol{\Sigma}_1) \|_2 \left(\| \boldsymbol{g}_i \|_2 + \frac{\| \boldsymbol{\mu} \|_2}{\sqrt{\lambda_1(\boldsymbol{\Sigma}_1)}} \right) - \frac{\| \boldsymbol{\mu} \|_2}{\sqrt{\lambda_1(\boldsymbol{\Sigma}_1)}}\Biggr\} \nonumber \\
& \qquad
\cap \Biggl\{ \|\boldsymbol{\gamma}_1(\boldsymbol{S}_{m}) - \boldsymbol{\gamma}_1(\boldsymbol{\Sigma}_1) \|_2 \leq \delta   \Biggr\} \biggl\arrowvert \mathrm{A}_i \Biggr) \nonumber  \\
& \quad + P\Biggl( \Biggl\{Z<  \| \boldsymbol{\gamma}_1(\boldsymbol{S}_{m}) - \boldsymbol{\gamma}_1(\boldsymbol{\Sigma}_1)\|_2 \left(\| \boldsymbol{g}_i \|_2 + \frac{\| \boldsymbol{\mu} \|_2}{\sqrt{\lambda_1(\boldsymbol{\Sigma}_1)}} \right) - \frac{\| \boldsymbol{\mu} \|_2}{\sqrt{\lambda_1(\boldsymbol{\Sigma}_1)}} \Biggr\} \nonumber\\
& \qquad \cap \Biggl\{ \|\boldsymbol{\gamma}_1(\boldsymbol{S}_{m}) - \boldsymbol{\gamma}_1(\boldsymbol{\Sigma}_1) \|_2 > \delta \Biggr\} \biggl\arrowvert \mathrm{A}_i \Biggr) \nonumber  \\ 
&\leq  P\Biggl( \Biggl\{Z<  \delta \left(\| \boldsymbol{g}_i \|_2 + \frac{\| \boldsymbol{\mu} \|_2}{\sqrt{\lambda_1(\boldsymbol{\Sigma}_1)}} \right) - \frac{\| \boldsymbol{\mu} \|_2}{\sqrt{\lambda_1(\boldsymbol{\Sigma}_1)}} \Biggr\} \nonumber \\
& \qquad \cap \Biggl\{ \|\boldsymbol{\gamma}_1(\boldsymbol{S}_{m}) - \boldsymbol{\gamma}_1(\boldsymbol{\Sigma}_1) \|_2 \leq \delta \Biggr\} \biggl\arrowvert \mathrm{A}_i \Biggr) 
+ P \Biggl( \|\boldsymbol{\gamma}_1(\boldsymbol{S}_{m}) - \boldsymbol{\gamma}_1(\boldsymbol{\Sigma}_1) \|_2 > \delta \biggl\arrowvert \mathrm{A}_i  \Biggr)  \nonumber  \\
& \leq
P\Biggl( Z<  \delta \left(\| \boldsymbol{g}_i \|_2 + \frac{\| \boldsymbol{\mu} \|_2}{\sqrt{\lambda_1(\boldsymbol{\Sigma}_1)}} \right) - \frac{\| \boldsymbol{\mu} \|_2}{\sqrt{\lambda_1(\boldsymbol{\Sigma}_1)}} \biggl\arrowvert \mathrm{A}_i  \Biggr) \nonumber \\
& \quad + P \Biggl( \|\boldsymbol{\gamma}_1(\boldsymbol{S}_{m}) - \boldsymbol{\gamma}_1(\boldsymbol{\Sigma}_1) \|_2 > \delta   \biggl\arrowvert \mathrm{A}_i \Biggr)\label{lem2p2}.
\end{align}
As for the first term on the right-hand side of \eqref{lem2p2}, we have
\begin{align*}
& P\Biggl( Z<  \delta \left(\| \boldsymbol{g}_i \|_2 + \frac{  \| \boldsymbol{\mu} \|_2}{ \sqrt{\lambda_1(\boldsymbol{\Sigma}_1)}} \right) - \frac{\| \boldsymbol{\mu} \|_2}{ \sqrt{\lambda_1(\boldsymbol{\Sigma}_1)} } \biggl\arrowvert \mathrm{A}_i \Biggr) \\
&=
P\Biggl( \Biggl\{Z<  \delta \left(\| \boldsymbol{g}_i \|_2 + \frac{  \| \boldsymbol{\mu} \|_2}{ \sqrt{\lambda_1(\boldsymbol{\Sigma}_1)}}\right) - \frac{\| \boldsymbol{\mu} \|_2}{ \sqrt{\lambda_1(\boldsymbol{\Sigma}_1)} } \Biggr\} \cap  \Biggl\{ \| \boldsymbol{g}_i \|_2 \leq c_1 \sqrt{n} \Biggr\} \biggl\arrowvert \mathrm{A}_i \Biggr) \\
&\quad + P\Biggl( \Biggl\{ Z<  \delta \left(\| \boldsymbol{g}_i \|_2 + \frac{  \| \boldsymbol{\mu} \|_2}{ \sqrt{\lambda_1(\boldsymbol{\Sigma}_1)}} \right) - \frac{\| \boldsymbol{\mu} \|_2}{ \sqrt{\lambda_1(\boldsymbol{\Sigma}_1)} } \Biggr\}  \cap  \Biggl\{ \| \boldsymbol{g}_i \|_2 > c_1 \sqrt{n}    \Biggr\} \biggl\arrowvert \mathrm{A}_i \Biggr) \\
&\leq P\Biggl( \Biggl\{Z<  \delta \left(\| \boldsymbol{g}_i \|_2 + \frac{  \| \boldsymbol{\mu} \|_2}{ \sqrt{\lambda_1(\boldsymbol{\Sigma}_1)}} \right) - \frac{\| \boldsymbol{\mu} \|_2}{ \sqrt{\lambda_1(\boldsymbol{\Sigma}_1)} } \Biggr\}  \cap  \Biggl\{ \| \boldsymbol{g}_i \|_2 \leq c_1 \sqrt{n}    \Biggr\} \biggl\arrowvert \mathrm{A}_i \Biggr) \\
& \quad+  P \Biggl( \| \boldsymbol{g}_i \|_2 > c_1 \sqrt{n} \biggl\arrowvert \mathrm{A}_i \Biggr) {.}
\end{align*}
The definition of $\delta$ yields
\begin{align*}
& P \Biggl( \Biggl\{Z<  \delta \left(\| \boldsymbol{g}_i \|_2 + \frac{  \| \boldsymbol{\mu} \|_2}{ \sqrt{\lambda_1(\boldsymbol{\Sigma}_1)}} \right) - \frac{\| \boldsymbol{\mu} \|_2}{\sqrt{\lambda_1(\boldsymbol{\Sigma}_1)}} \Biggr\} \cap \Biggl\{\|\boldsymbol{g}_i \|_2 \leq c_1 \sqrt{n} \Biggr\} \biggl\arrowvert \mathrm{A}_i \Biggr)  \\
&\leq 
P \Biggl( \Biggl\{Z<  \delta \left(c_1 \sqrt{n} + \frac{  \| \boldsymbol{\mu} \|_2}{ \sqrt{\lambda_1(\boldsymbol{\Sigma}_1)}} \right) - \frac{\| \boldsymbol{\mu} \|_2}{\sqrt{\lambda_1(\boldsymbol{\Sigma}_1)}} \Biggr\} \cap \Biggl\{\|\boldsymbol{g}_i \|_2 \leq c_1 \sqrt{n} \Biggr\} \biggl\arrowvert \mathrm{A}_i \Biggr)  
\\
&\leq
P \Biggl( Z<  \delta \left(c_1 \sqrt{n} + \frac{  \| \boldsymbol{\mu} \|_2}{ \sqrt{\lambda_1(\boldsymbol{\Sigma}_1)}} \right) - \frac{\| \boldsymbol{\mu} \|_2}{\sqrt{\lambda_1(\boldsymbol{\Sigma}_1)}} \biggl\arrowvert \mathrm{A}_i \Biggr)  
\\
&\leq P \Biggl(Z< \frac{-(1-\alpha) \| \boldsymbol{\mu} \|_2}{\sqrt{\lambda_1(\boldsymbol{\Sigma}_1)}} \biggl\arrowvert \mathrm{A}_i \Biggr) .
\end{align*}
Moreover, recalling that
$c_1  = 1 + { K_g^2}/{\sqrt{c}}$,
we have
\begin{align*}
&P(\|\boldsymbol{g}_i\|_2 > c_1 \sqrt{n}  \arrowvert \mathrm{A}_i ) 
= P( |  \| \boldsymbol{g}_i \|_2 - \sqrt{n}  + \sqrt{n} | > c_1 \sqrt{n} \arrowvert \mathrm{A}_i ) \\
& \leq P( |  \| \boldsymbol{g}_i \|_2 - \sqrt{n}  | > (c_1-1) \sqrt{n} \arrowvert \mathrm{A}_i ) \\
& \leq 2 \e^{-n}, 
\end{align*}
where \eqref{v33} is used in the last inequality.
It follows that
\begin{align*}
&P\Biggl( Z<  \delta \left(\| \boldsymbol{g}_i \|_2 + \frac{  \| \boldsymbol{\mu} \|_2}{ \sqrt{\lambda_1(\boldsymbol{\Sigma}_1)}} \right) - \frac{\| \boldsymbol{\mu} \|_2}{ \sqrt{\lambda_1(\boldsymbol{\Sigma}_1)} } \biggl\arrowvert \mathrm{A}_i \Biggr) \\ 
& \leq P\Biggl(Z< \frac{-(1-\alpha) \| \boldsymbol{\mu} \|_2}{\sqrt{\lambda_1(\boldsymbol{\Sigma}_1)}} \biggl\arrowvert \mathrm{A}_i \Biggr) +  2 \e^{-n} .
\end{align*}
As for the second term on the right-hand side of \eqref{lem2p2}, by using \eqref{AE1}, \eqref{AE2}, \eqref{AE3} and the Davis--Kahan theorem (see, e.g., \citet[p.89]{RefV18}), we have
\begin{align*}
& P\left(\|\boldsymbol{\gamma}_1(\boldsymbol{S}_{m}) - \boldsymbol{\gamma}_1(\boldsymbol{\Sigma}_1) \|_2 > \delta  \biggl\arrowvert \mathrm{A}_i  \right)
= 
P\left(\|\boldsymbol{\gamma}_1(\boldsymbol{S}_{m}) - \boldsymbol{\gamma}_1(\boldsymbol{\Sigma}) \|_2 > \delta  \biggl\arrowvert \mathrm{A}_i  \right) \\
&\leq 
 P\left( \frac{2^{{3}/2}\| \boldsymbol{S}_{m} - \boldsymbol{\Sigma}\|_{\mathrm{op}}}{\lambda_1(\boldsymbol{\Sigma})- \lambda_2(\boldsymbol{\Sigma})}> \delta \biggl\arrowvert \mathrm{A}_i \right) \\
&\leq
P\left( \frac{2^{{3}/2}\| \boldsymbol{S}_{m} - \boldsymbol{\Sigma}\|_{\mathrm{op}}}{\|\boldsymbol{\mu}\|_2^2 }> \delta \biggl\arrowvert \mathrm{A}_i \right) =
 P \left( \| \boldsymbol{S}_{m} - \boldsymbol{\Sigma} \|_{\mathrm{op}} > \frac{\delta \|\boldsymbol{\mu} \|_2^2}{2^{{3}/{2}}} \biggl\arrowvert \mathrm{A}_i \right).
\end{align*}
The definition of $\delta$ and Proposition~\ref{evs} with $u = n$ yield
   \begin{align*} 
& P \left( \| \boldsymbol{S}_{m} - \boldsymbol{\Sigma} \|_{\mathrm{op}} > \frac{\delta \|\boldsymbol{\mu} \|_2^2}{2^{{3}/{2}}} \biggl\arrowvert \mathrm{A}_i \right)\\*
&\leq
 P \biggr( \| \boldsymbol{S}_{m} - \boldsymbol{\Sigma} \|_{\mathrm{op}} > C K^2 \left(\sqrt{\frac{2n}{m}}+\frac{2n}{m} \right)  \left( \frac{\lambda_1(\boldsymbol{\Sigma}_1)}{2} + \frac{\lambda_1(\boldsymbol{\Sigma}_2)}{2}  +  \| \boldsymbol{\mu} \|_2^2  \right) \biggl\arrowvert \mathrm{A}_i \biggr) \\*
&=
\frac{1}{P \left( \mathrm{A}_i \right)} \\
& \quad \cdot P \left( \| \boldsymbol{S}_{m} - \boldsymbol{\Sigma} \|_{\mathrm{op}} > C K^2 \left(\sqrt{\frac{2n}{m}}+\frac{2n}{m} \right) \left( \frac{\lambda_1(\boldsymbol{\Sigma}_1)}{2} + \frac{\lambda_1(\boldsymbol{\Sigma}_2)}{2}  +  \| \boldsymbol{\mu} \|_2^2  \right) \cap \mathrm{A}_i \right) \\*
&\leq
 2 P \biggr( \| \boldsymbol{S}_{m} - \boldsymbol{\Sigma} \|_{\mathrm{op}} > C K^2 \left(\sqrt{\frac{2n}{m}}+\frac{2n}{m} \right) \left( \frac{\lambda_1(\boldsymbol{\Sigma}_1)}{2} + \frac{\lambda_1(\boldsymbol{\Sigma}_2)}{2}  +  \| \boldsymbol{\mu} \|_2^2  \right) \biggr)\\
&
\leq 4 \e^{-n} .
\end{align*}
From what has already been proved, we conclude that
\begin{equation*}
P \biggl(\langle \boldsymbol{\gamma}_1(\boldsymbol{S}_{m}), \boldsymbol{X}_i\rangle<0 \biggl\arrowvert \mathrm{A}_i \biggr) \\
\leq P\biggl(Z< \frac{-(1-\alpha) \| \boldsymbol{\mu} \|_2}{\sqrt{\lambda_1(\boldsymbol{\Sigma}_1)}} \biggl\arrowvert \mathrm{A}_i \biggr) +  6 \e^{-n}.
\end{equation*}
This completes the proof.
\qed

\subsection{Proof of Corollary \ref{col}}
Consider $m$, $n$, and $\eta$ satisfying \eqref{thas}.
Then, by using the Bonferroni inequality and \eqref{millsa}, we have
\begin{align}
&P \left( \left\{ \bigcap_{i=1}^m \bigl\{ \theta_i \langle \boldsymbol{\gamma}_1(\boldsymbol{S}_{m}), \boldsymbol{X}_i \rangle > 0 \bigr\} \right\} \cup \left\{ \bigcap_{i=1}^m \bigl\{ \theta_i \langle \boldsymbol{\gamma}_1(\boldsymbol{S}_{m}), \boldsymbol{X}_i \rangle <0 \bigr\} \right\}  \right) 
\nonumber \\
&\geq
P \left(\bigcap_{i=1}^m \bigl\{ \theta_i \langle \boldsymbol{\gamma}_1(\boldsymbol{S}_{m}), \boldsymbol{X}_i \rangle > 0 \bigr\} \right) 
= 1 - P \left(\bigcup_{i=1}^m \bigl\{ \theta_i \langle \boldsymbol{\gamma}_1(\boldsymbol{S}_{m}), \boldsymbol{X}_i \rangle \leq 0 \bigr\} \right) 
\nonumber \\
&\geq
1 - \sum_{i=1}^m P \left( \bigl\{ \theta_i \langle \boldsymbol{\gamma}_1(\boldsymbol{S}_{m}), \boldsymbol{X}_i \rangle \leq 0 \bigr\} \right)
\nonumber \\
&\geq
1 - \frac{m}{\sqrt{2\pi (1-\alpha)^{2}\eta}} \exp\left(- \frac{(1-\alpha)^2 \eta}{2} \right) - 6 m \e^{-n}.
\label{col:eq}
\end{align}
The second and third terms on the right-hand side of \eqref{col:eq} converge to 0 as $n,m \to\infty$ with \eqref{ARh} under the assumption that $n/\eta=O(1)$. 
\qed

\section{Concluding remarks}\label{sec:5}
In this paper, we derived non-asymptotic bounds for the error probability of the spectral clustering algorithm when the mixture distribution of two multivariate normal distributions that form the allometric extension relationship is considered.
As future directions, it is interesting to relax the assumption of the normal distribution to the sub-gaussian distribution and to consider weights in the mixture distribution other than $\pi_1=\pi_2=1/2$.

\begin{acks}[Acknowledgments]
This study was supported in part by Japan Society for the Promotion of Science KAKENHI Grant Numbers 21K13836 and 23K16851.
\end{acks}


\begin{thebibliography}{99}	
\bibitem[Abbe et al.(2022)]{RefAFW22}
Abbe, E., Fan, J., Wang, K. (2022).
An $\ell_p$ theory of PCA and spectral clustering.
\textit{Ann. Statist.} \textbf{50}, no.4, 2359--2385.
%
\bibitem[Amit et al.(2017)]{RefA17}
Amit, S., Mukesh, P., Akshansh, G., Neha, B., Om, P.~P., Aruna, T., Meng, J.~E., Weiping, D., Chin-Teng, L. (2017).
A review of clustering techniques and developments.
\textit{Neurocomputing} \textbf{267}, 664--681.
%
\bibitem[Bartoletti et al.(1999) ]{RefBFN99} 
Bartoletti, S., Flury, B.~D., Nel, D.~G. (1999).
Allometric extension.
\textit{Biometrics} \textbf{55}, no.4, 1210--1214.
%
\bibitem[Borysov et al.(2014)]{RefB14}
Borysov, P., Hannig, J., Marron, J. S. (2014).
Asymptotics of hierarchical clustering for growing dimension.
\textit{J. Multivariate Anal.} \textbf{124}, 465--479.
%
\bibitem[Cai and Zhang(2018)]{RefCZ18}
Cai, T.~T., Zhang, A. (2018).
Rate-optimal perturbation bounds for singular subspaces with applications to high-dimensional statistics.
\textit{Ann. Statist.} \textbf{46}, no.1, 60--89.
%
\bibitem[Flury(1997)]{RefF97} 
Flury, B. (1997).
\textit{A First Course in Multivariate Statistics}.
Springer-Verlag, New York.
%
\bibitem[Hills(2006)]{RefH06} 
Hills, M. (2006). 
Allometry. In \textit{Encyclopedia of Statistical Sciences} (eds S. Kotz, C.B. Read, N. Balakrishnan, B. Vidakovic and N.L. Johnson). 
https://doi.org/10.1002/0471667196.ess0033.pub2
(Last access: 2023/06/01)
%
\bibitem[Hsu and Kakade(2013)]{RefDS13}
Hsu, D., Kakade, S. (2013).
Learning mixtures of spherical Gaussians: moment methods and spectral decompositions,
in \textit{ITCS'13---Proceedings of the 2013 ACM Conference on Innovations in Theoretical Computer Science}. ACM, New York. 
pp.11-19.
%
\bibitem[Kurata et al.(2008)]{RefKHF08}
Kurata, H., Hoshino, T., Fujikoshi, Y. (2008).
Allometric extension model for conditional distributions.
\textit{J. Multivariate Anal.} \textbf{99}, no.9, 1985--1998.
%
\bibitem[L\"{o}ffler et al.(2021)]{RefL21}
L\"{o}ffler, M., Zhang, A.~Y., Zhou, H.~H. (2021).
Optimality of spectral clustering in the Gaussian mixture model
\textit{Ann. Statist.} \textbf{49}, no.5, 2506--2530.
%
\bibitem[Matsuura and Kurata(2014)]{RefMK14}
Matsuura, S., Kurata, H. (2014).
Principal points for an allometric extension model.
\textit{Statist. Papers} \textbf{55}, no.3, 853--870.
%
\bibitem[Ndaoud(2022)]{RefN22}
Ndaoud, M. (2022).
Sharp optimal recovery in the two component Gaussian mixture model.
\textit{Ann. Statist.} \textbf{50}, no.4, 2096--2126.
%
\bibitem[O'Neill(1978)]{RefO78} 
O'Neill, T.~J. (1978).
Normal discrimination with unclassified observations.
\textit{J. Amer. Statist. Assoc.} \textbf{73}, no.364, 821--826.
%
\bibitem[Pollard(1981)]{RefP81}
Pollard, D. (1981).
Strong consistency of $k$-means clustering.
\textit{Ann. Statist.} \textbf{9}, no.1, 135--140.
%
\bibitem[Pollard(1982)]{RefP82}
Pollard, D. (1982).
A central limit theorem for $k$-means clustering.
\textit{Ann. Probab.} \textbf{10}, no.4, 919--926.
%
\bibitem[Tarpey(2007)]{RefT07}
Tarpey, T. (2007).
Linear transformations and the $k$-means clustering algorithm.
\textit{Amer. Statist.} \textbf{61}, no.1, 34--40.
%
\bibitem[Tsukuda and Matsuura(2023)]{RefTM23}
Tsukuda, K., Matsuura, S. (2023).
High-dimensional hypothesis testing for allometric extension model.
\textit{J. Multivariate Anal.} \textbf{197}, 105208.
%
\bibitem[Vershynin(2018)]{RefV18} 
Vershynin, R. (2018).
\textit{High-Dimensional Probability. An Introduction with Applications in Data Science}. 
Cambridge University Press, Cambridge.
%
\end{thebibliography}
\end{document}